# DECISION THEORY RESULTS FOR ONE-SIDED MULTIPLE COMPARISON PROCEDURES[1]

BY ARTHUR COHEN AND HAROLD B. SACKROWITZ

*Rutgers University*

A resurgence of interest in multiple hypothesis testing has occurred in the last decade. Motivated by studies in genomics, microarrays, DNA sequencing, drug screening, clinical trials, bioassays, education and psychology, statisticians have been devoting considerable research energy in an effort to properly analyze multiple endpoint data. In response to new applications, new criteria and new methodology, many ad hoc procedures have emerged. The classical requirement has been to use procedures which control the strong familywise error rate (FWE) at some predetermined level $\alpha$. That is, the probability of any false rejection of a true null hypothesis should be less than or equal to $\alpha$. Finding desirable and powerful multiple test procedures is difficult under this requirement.

One of the more recent ideas is concerned with controlling the false discovery rate (FDR), that is, the expected proportion of rejected hypotheses which are, in fact, true. Many multiple test procedures do control the FDR.

A much earlier approach to multiple testing was formulated by Lehmann [*Ann. Math. Statist.* **23** (1952) 541–552 and **28** (1957) 1–25]. Lehmann's approach is decision theoretic and he treats the multiple endpoints problem as a $2^k$ finite action problem when there are $k$ endpoints. This approach is appealing since unlike the FWE and FDR criteria, the finite action approach pays attention to false acceptances as well as false rejections. In this paper we view the multiple endpoints problem as a $2^k$ finite action problem. We study the popular procedures single-step, step-down and step-up from the point of view of admissibility, Bayes and limit of Bayes properties. For our model, which is a prototypical one, and our loss function, we are able to demonstrate the following results under some fairly general conditions to be specified:

Received February 2003; revised March 2004.
[1]Supported by NSA Grant MDA 904-02-1-0039.
*AMS 2000 subject classifications.* 62F03, 62C15.
*Key words and phrases.* Global problem, single-step procedure, step-down procedure, step-up procedure, Bayes procedure, admissibility, familywise error rate, false acceptance rate, false discovery rate, partitioning principle, closure property, finite action problem.

(i) The single-step procedure is admissible.

(ii) A sequence of prior distributions is given for which the step-down procedure is a limit of a sequence of Bayes procedures.

(iii) For a vector risk function, where each component is the risk for an individual testing problem, various admissibility and inadmissibility results are obtained.

In a companion paper [Cohen and Sackrowitz, *Ann. Statist.* **33** (2005) 145–158], we are able to give a characterization of Bayes procedures and their limits. The characterization yields a complete class and the additional useful result that the step-up procedure is inadmissible. The inadmissibility of step-up is demonstrated there for a more stringent loss function. Additional decision theoretic type results are also obtained in this paper.

**1. Introduction.** A resurgence of interest in multiple hypothesis testing has occurred in the last decade. Motivated by studies in genomics, microarrays, DNA sequencing, drug screening, clinical trials, bioassays, education and psychology, statisticians have been devoting considerable research energy in an effort to properly analyze multiple endpoint data. In response to new applications, new criteria and new methodology, many ad hoc procedures have emerged. The classical requirement has been to use procedures which control the strong familywise error rate (FWE) at some predetermined level $\alpha$. That is, the probability of any false rejection of a true null hypothesis should be less than or equal to $\alpha$. Finding desirable and powerful multiple test procedures is difficult under this requirement. Two useful tools for the construction of such multiple level-$\alpha$ tests are the closure principle [see Marcus, Peritz and Gabriel (1976), as well as Hochberg and Tamhane (1987)] and the partitioning principle [see Stefánsson, Kim and Hsu (1988) and Finner and Strassburger (2002)]. These tools can be used to generate large classes of multiple test procedures satisfying the FWE criterion.

One of the more recent ideas is concerned with controlling the false discovery rate (FDR), that is, the expected proportion of rejected hypotheses which are, in fact, true. Many multiple test procedures do control the FDR. See, for example, Benjamini and Hochberg (1995), Benjamini and Yekutieli (2001), Efron, Tibshirani, Storey and Tusher (2001) and Sarkar (2002). This criterion is particularly appealing if the number of endpoints is large. In some modern applications this number can be in the thousands. A summary of studies on multiple endpoint methods used with microarray data is given in Dudoit, Shaffer and Boldrick (DSB) (2003).

The closure and partitioning principles tend to be linked to the step-down approach described in Hochberg and Tamhane (1987) and studied extensively in the literature. FDR was initially linked to the step-up approach. See Hochberg (1988). More recently, step-down and combined step-down with step-up methods have been viewed from an FDR point of view. Sarkar (2002) notes: "While the FDR has been receiving increasing attention by researchers in different fields of statistics, theoretical progress has not been made at a similar pace." Sarkar's remark applies to the entire area of multiple endpoint testing. Finner and Strassburger (2002) say: "Further and in general difficult problems are the comparison of different multiple test procedures and the related questions concerning admissibility." They go on to say: "A serious issue is optimality and admissibility of multiple decision procedures." DSB (2003) remark "Optimality of multiple tests is an interesting research avenue to pursue from both a theoretical and a practical point of view."

A much earlier approach to multiple testing was formulated by Lehmann (1952, 1957). Lehmann's approach is decision theoretic and he treats the multiple endpoints problem as a $2^k$ finite action problem when there are $k$ endpoints. The formulation as a $2^k$ action problem is particularly appealing since in terms of what is desired, one wishes to decide whether to accept or reject for each of the $k$ hypotheses posed. This approach entails the specification of losses, which can be quite general. Lehmann (1952, 1957) was able to demonstrate some optimality properties for the single-step procedure and step-down procedure in two-dimensional problems for some hypotheses and for some restricted classes of procedures. Methods developed through the years to further the theory of testing a single hypothesis (a two-action problem) do not extend easily to multiple actions and little progress has been made for this model. Nevertheless, the potential and importance of this approach are compelling since the evaluation of methodologies and procedures is wanting and necessary in this subject. Little is known about properties of the various procedures and rigorously studying the underpinnings of the methodologies is essential. Furthermore, unlike the FWE and FDR criteria, the finite action approach pays attention to false acceptances as well as false rejections.

Our approach will be to regard the problem as a $2^k$ action problem. We carefully distinguish between what is known as the global problem and multiple endpoints problem. We are very precise about what null hypotheses and what alternative hypotheses are to be considered. Several notions of monotonicity of procedures and monotonicity of risk functions have been introduced and studied in Cohen and Sackrowitz (CS) (2004). In this paper Bayes procedures, limits of Bayes procedures and admissibility results of procedures are studied. In particular we examine single-step, step-down and step-up procedures. We note that DSB (2003) classify the 18 procedures they study as single step, step-down or step-up. We consider loss functions that are sums of losses for the individual endpoints.

We confine our study to a simple but prototypical model, although many of the results would remain true for other models. The model assumed is that we observe a $(k \times 1)$ random vector $\mathbf{Z}$ which is assumed to be $k$-variate normal with mean vector $\boldsymbol{\mu}$ and known covariance $\Sigma$. Among the results are the following:

RESULT 1.1. Suppose the covariance matrix $\Sigma$ is of the intraclass type, that is, all variances equal, all covariances equal. Then under some mild conditions the single-step procedure is admissible. The approach used to prove admissibility is somewhat new.

RESULT 1.2. If $\Sigma = \sigma^2 I$, the step-down procedure studied is shown to be a limit of a sequence of Bayes procedures.

RESULT 1.3. Suppose $\Sigma$ is intraclass, and $\rho$ is the common correlation coefficient between any pair of variables. Consider a vector risk (VRI) where the components of the vector are the risks for the individual testing problems. Then the single-step procedure is admissible for any $-1 < \rho < 1$. The step-down and step-up procedures are admissible if and only if $\rho \geq 0$. As a corollary it follows that for $\rho < 0$, step-up and step-down are inadmissible for the loss function which sums the losses for the individual component problems.

In Section 2 we state the models, distinguish between global test problems and multiple endpoint testing problems, introduce the loss functions, describe the various properties of procedures and give other preliminaries. In Section 3 we describe the single-step, step-down and step-up procedures. In Sections 4 and 5 we state properties of these procedures. All proofs are given in the Appendix.

## 2. Models and preliminaries.

2.1. *Models.* Let $\mathbf{Z}$ be a $(k \times 1)$ random vector which is $k$-variate normal with mean vector $\boldsymbol{\mu}$ and known covariance $\Sigma$. One global one-sided hypothesis testing problem is

$$(2.1) \qquad H^{(G)} : \boldsymbol{\mu} = \mathbf{0} \quad \text{vs} \quad K^{(G)} : \boldsymbol{\mu} \geq \mathbf{0} \setminus \{\boldsymbol{\mu} = \mathbf{0}\},$$

that is, $\mu_i \geq 0$, $i = 1, 2, \ldots, k$, with at least one $\mu_i > 0$. Such a problem is distinguished from a one-sided multiple endpoints problem in which one tests

$$(2.2) \qquad H_i : \mu_i = 0 \quad \text{vs} \quad K_i : \mu_i > 0, \qquad i = 1, 2, \ldots, k.$$

That is, the latter problem is a $2^k$ action problem where one selects an action to either accept or reject $H_i$, $i = 1, 2, \ldots, k$.

Note that another form of the one-sided multiple endpoints problem is

$$(2.3) \qquad H_i^* : \mu_i \leq 0 \quad \text{vs} \quad K_i : \mu_i > 0.$$

In the multiple endpoints literature there are ample instances of both the $H_i$ and $H_i^*$ problems. We mention both since from a decision theory point of view we will see that sometimes different results ensue depending on whether $H_i$ or $H_i^*$ is being tested. In connection with distinguishing between $H_i$ and $H_i^*$ we mention two practical situations where the multiple endpoints scenario arises.

(I) Consider the problem of comparing $k$ treatments with a control assuming (often realistically) that the treatment mean will be at least as large as the control mean. [This model is called the tree order model in Robertson, Wright and Dykstra (1988).] Then if $Z_i$ represents the difference between a reading on the $i$th treatment and the control, $Z_i$ has mean $\mu_i$ where $\mu_i \geq 0$. Assuming all treatment and control observations are normal, independent, with variances 1, then $\mathbf{Z}$ is multivariate normal with mean vector $\boldsymbol{\mu}$ and covariance matrix $\Sigma$. The covariance matrix $\Sigma$ is $(k \times k)$ and

$$(2.4) \qquad \Sigma = 2 \begin{pmatrix} 1 & 1/2 & \cdots & 1/2 \\ 1/2 & 1 & \cdots & 1/2 \\ \vdots & & & \vdots \\ 1/2 & & & 1 \end{pmatrix}.$$

The appropriate multiple hypotheses in this case are those in (2.2).

Note that the $(k \times k)$ covariance matrix $\Sigma$ in (2.4) is a special case of a class of covariance matrices which are called intraclass. That is, a covariance matrix $\Sigma = (\sigma_{ij})$ is intraclass if $\sigma_{ii}$ are the same for all $i = 1, 2, \ldots, k$, and $\sigma_{ij}$, $i = 1, \ldots, k, j = 1, \ldots, k, i \neq j$, are the same. When $\Sigma$ is intraclass, then the normal variables are exchangeable. Note also that a special case of intraclass is when all $\sigma_{ii}$ are the same and all $\sigma_{ij} = 0$, $i \neq j$. In this latter case the $Z_i$ are independent. Another special case of intraclass is when $k = 2$, and $\sigma_{11} = \sigma_{22}$. The intraclass matrix $\Sigma$ may be written as

$$(2.5) \qquad \Sigma = \sigma^2 \begin{pmatrix} 1 & \rho & \cdots & \rho \\ \rho & 1 & \cdots & \rho \\ \vdots & & & \vdots \\ \rho & \cdots & \cdots & 1 \end{pmatrix}$$

with $\rho$ restricted to the interval $-1/(k-1) \leq \rho \leq 1$. See, for example, Krishnaiah and Pathak (1967).

(II) Let $\mathbf{X}_i$, $i = 1, 2$, be independent normal with mean vector $\boldsymbol{\nu}_i$ and known covariance $\Sigma_i$. $\mathbf{X}_1$ corresponds to a $(k \times 1)$ vector of measurements made on a control subject. $\mathbf{X}_2$ corresponds to a $(k \times 1)$ vector of measurements made on a treatment subject. Consider $\mathbf{Z} = \mathbf{X}_2 - \mathbf{X}_1$ and note $\mathbf{Z}$ is multivariate normal with mean vector $\boldsymbol{\mu} = \boldsymbol{\nu}_2 - \boldsymbol{\nu}_1$ and covariance matrix $\Sigma = \Sigma_1 + \Sigma_2$. If one feels that the treatment cannot decrease $\nu_{1i}$, $i = 1, 2, \ldots, k$, then this is the classic multiple endpoints problem with (2.2) as the multiple hypotheses. If one feels that the treatment can reduce $\nu_{1i}$ as well as increase $\nu_{1i}$, then this is the classic multiple endpoints problem with (2.3) as the hypotheses.

2.2. *Preliminaries.* A $2^k$ finite action problem has actions $\mathbf{a} = (a_1, a_2, \ldots, a_k)'$ where $a_i$ equals 0 or 1 for $i = 1, \ldots, k$. An action where $a_i = 1$ means that $H_i$ is rejected, where if $a_i = 0$, $H_i$ is accepted. Thus, for example, $\mathbf{a} = (1, \ldots, 1)'$ means all $H_i$ are rejected. It will be convenient to define

$$\Gamma = \{\mathbf{u} : \mathbf{u} = (u_1, \ldots, u_k)', u_i = 0 \text{ or } 1, \text{ all } i\}.$$

Note that $\Gamma$ can be used to represent the totality of all actions. However, $\Gamma$ will serve other purposes as well.

Decision rules $\delta(\cdot|\mathbf{z})$ are probability mass functions on $\Gamma$ with the interpretation that $\delta(\mathbf{a}|\mathbf{z})$ is the conditional probability of action $\mathbf{a}$ given $\mathbf{z}$ is observed. For each $\mathbf{z}$, a nonrandomized decision rule chooses a single element of $\Gamma$ with probability 1 and assigns all other actions probability 0. Each decision rule $\delta$ determines a set of test functions for the individual testing problems. These test functions are given by $\boldsymbol{\psi}(\mathbf{z}) = (\psi_1(\mathbf{z}), \ldots, \psi_k(\mathbf{z}))'$ where $\psi_i(\mathbf{z})$ is the probability of rejecting $H_i$. A decision procedure $\delta(\mathbf{a}|\mathbf{z})$ determines a set of $\psi_i^\delta(\mathbf{z})$, $i = 1, \ldots, k$, as follows:

$$(2.6) \qquad \psi_i^\delta(\mathbf{z}) = \sum_{\mathbf{a} \in A_i} \delta(\mathbf{a}|\mathbf{z}) = \sum_{\mathbf{a} \in \Gamma} a_i \delta(\mathbf{a}|\mathbf{z}),$$

where $A_i = \{\mathbf{a} \in \Gamma : \mathbf{a} \text{ has a 1 in the } i\text{th position}\}$. Whereas $\delta(\mathbf{a}|\mathbf{z})$ deter- mines $\boldsymbol{\psi}(\mathbf{z})$, the reverse is not true. If $\boldsymbol{\psi}(\mathbf{z})$ is nonrandomized it uniquely determines some $\delta(\mathbf{a}|\mathbf{z})$. The $\delta(\mathbf{a}|\mathbf{z})$ determined is nonrandomized.

For problem (2.2), we partition the parameter space $\Omega = \{\boldsymbol{\mu} : \mu_i \geq 0, i = 1, \ldots, k\}$ into $2^k$ sets $\Omega_\mathbf{v}$, $\mathbf{v} \in \Gamma$, where $\Omega_\mathbf{v} = \{\boldsymbol{\mu} : \boldsymbol{\mu} = (\mu_1, \mu_2, \ldots, \mu_k)'$, and $\mu_i > 0$ if $v_i = 1$ and $\mu_i = 0$ if $v_i = 0$, $i = 1, \ldots, k\}$. For problem (2.3), $\Omega = \mathbb{R}^k$ and we have a similar partition but $v_i = 0$ means $\mu_i \leq 0$. Also for problem (2.2) let $\Omega^{(i)} = \{\boldsymbol{\mu} : \boldsymbol{\mu} \in \Omega, \mu_i = 0\}$.

A loss function is a function of the action taken and the true state of nature. We will take the loss function to be additive over the individual component problems and for each component problem we choose the loss as follows: zero loss for a correct decision; a loss of 1 for rejecting $H_i$ when it is true and a loss of $b$ for accepting $H_i$ when it is false. The loss function for the finite action problem can be expressed as

$$(2.7) \qquad L(\mathbf{a}, \boldsymbol{\mu}) = \sum_{i=1}^k a_i(1 - v_i) + \sum_{i=1}^k b(1 - a_i)v_i, \qquad \boldsymbol{\mu} \in \Omega_\mathbf{v},$$

with $0 < b$. This loss function reflects the property that the loss is additive over the losses for the component problems.

The risk function for a decision procedure $\delta$ is

(2.8) $$R(\delta, \boldsymbol{\mu}) = E_{\boldsymbol{\mu}} \sum_{\mathbf{a} \in \Gamma} L(\mathbf{a}, \boldsymbol{\mu}) \delta(\mathbf{a}|\mathbf{z}).$$

For the above loss function (2.7) it follows from (2.6) that the risk depends on $\delta$ only through $\boldsymbol{\psi}$, so we can write (2.8) as

(2.9) $$R(\boldsymbol{\psi}, \boldsymbol{\mu}) = \sum_{i=1}^{k} R_{(i)}(\psi_i, \boldsymbol{\mu}),$$

where

(2.10) $$R_{(i)}(\psi_i, \boldsymbol{\mu}) = \begin{cases} E_{\boldsymbol{\mu}} \psi_i(\mathbf{z}), & \mu_i = 0, \\ b(1 - E_{\boldsymbol{\mu}} \psi_i(\mathbf{z})), & \mu_i > 0. \end{cases}$$

The risk function (2.9) can be written as

(2.11) $$E_{\boldsymbol{\mu}}(\boldsymbol{\psi}'(\mathbf{1} - \mathbf{v}) + b(\mathbf{1} - \boldsymbol{\psi})' \mathbf{v}),$$

where $\mathbf{1} = (1, \ldots, 1)'$.

A decision procedure $\boldsymbol{\psi}$ is said to be inadmissible if there exists another procedure $\boldsymbol{\psi}^*$ such that $R(\boldsymbol{\psi}^*, \boldsymbol{\mu}) \leq R(\boldsymbol{\psi}, \boldsymbol{\mu})$ for every $\boldsymbol{\mu}$ with strict inequality for some $\boldsymbol{\mu}$. Otherwise $\boldsymbol{\psi}$ is admissible.

As previously noted we can view the multiple endpoints problem as one involving $k$ endpoints in which $\psi_i(\mathbf{z})$, $i = 1, \ldots, k$, is a test function for the $i$th endpoint. In this scenario one may wish to consider a vector risk approach where the risk consists of a $(k \times 1)$ vector $\mathbf{R} = (R_{(1)}, \ldots, R_{(k)})'$, $R_{(i)} = R(\psi_i, \boldsymbol{\mu})$ given in (2.10). In this formulation any procedure which has an admissible test for each single component is admissible in the vector risk formulation. For general results concerning vector risks, see Cohen and Sackrowitz (1984).

**3. Some procedures for multiple endpoint problems.** We focus on three special cases of the most frequently discussed procedures, namely single-step, step-down and step-up. The three procedures are considered in Hochberg and Tamhane (1987) and Shaffer (1995). In all that follows we assume without loss of generality that the variance of each $Z_i$ is 1. Furthermore, for now for step-up and step-down we limit our discussion to procedures which are symmetric in the variables $Z_1, \ldots, Z_k$, that is, procedures that are permutation equivariant. The normal model, with intraclass covariance matrix, represents the most general case of permutation invariance.

3.1. *Single-step.* The single-step procedure we study is:

PROCEDURE 3.1. Reject $H_i$ if and only if $Z_i > C_i$.

The constants $C_i$ are typically chosen so that the strong familywise error rate (FWE) is less than or equal to $\alpha$.

3.2. *Step-down.* The step-down procedure we study is as follows:

PROCEDURE 3.2. Let $Z_{(1)} \leq Z_{(2)} \leq \cdots \leq Z_{(k)}$ be the order statistics for the set $(Z_1, Z_2, \ldots, Z_k)'$ and let $C_j$ be a strictly increasing set of critical values:

(i) If $Z_{(k)} > C_k$, reject $H_{(k)}$. Otherwise accept all $H_i$.
(ii) If $H_{(k)}$ is rejected, reject $H_{(k-1)}$ if $Z_{(k-1)} > C_{k-1}$. Otherwise accept all $H_{(k-1)}, \ldots, H_{(1)}$.
(iii) In general, at stage $j$, if $Z_{(j)} > C_j$, reject $H_{(j)}$. Otherwise accept $H_{(j)}, \ldots, H_{(1)}$.

The critical values may be chosen so that:

OUTCOME 3.3. $P\{Z_{(k)} \leq C_k\} = 1 - \alpha$ when all $H_i$ are true; $P\{Z_{(k-1)} \leq C_{k-1}\} = 1 - \alpha$, with $Z_{(k)}$ excluded, and $H_{(1)}, \ldots, H_{(k-1)}$ are true. That is, after one of the hypotheses is rejected, we consider a new problem with the $(k-1)$ remaining variables that correspond to those parameters not rejected at step 1. $P\{Z_{(j)} \leq C_j\} = 1 - \alpha$, when $(k - j)$ variables and their corresponding hypotheses are excluded and $H_{(1)}, \ldots, H_{(j)}$ are true. This choice of constants leads to control of the strong FWE.

The step-down procedure results by applying the closure method. See Hochberg and Tamhane [(1987), Chapter 2, Section 4.1] for a description of this method. The method is used by Marcus, Peritz and Gabriel (1976). Finner and Roters (2002) note that this method strongly controls the FWE and they call such a testing method a multiple level-$\alpha$ test procedure. The step-down method rejects more $H_i$'s than the single-step procedure for the same $\alpha$. The single-step procedure would use $C_k$ for $C_j$, all $j = 1, \ldots, k$.

REMARK. Procedure 3.2 is one type of step-down procedure. Another type, used, for example, by Marcus, Peritz and Gabriel (1976), uses a likelihood ratio test in applying the closure method. This results in a different step-down procedure.

3.3. *Step-up.* The step-up procedure we study is as follows:

PROCEDURE 3.4. Let $Z_{(1)} \leq Z_{(2)} \leq \cdots \leq Z_{(k)}$ be the order statistics for the set $(Z_1, \ldots, Z_k)$ and let $C_j$ be a strictly increasing set of critical values.

(i) If $Z_{(1)} \leq C_1$, accept $H_{(1)}$. Otherwise reject all $H_i$.
(ii) If $H_{(1)}$ is accepted, accept $H_{(2)}$ if $Z_{(2)} \leq C_2$. Otherwise reject $H_{(2)}, \ldots, H_{(k)}$.
(iii) In general, at stage $j$, if $Z_{(j)} \leq C_j$, accept $H_{(j)}$. Otherwise reject $H_{(j)}, \ldots, H_{(k)}$.

The critical values $C_j$ are sometimes chosen so that:

OUTCOME 3.5. $P\{Z_{(1)} \leq C_1, Z_{(2)} \leq C_2, \ldots, Z_{(j)} \leq C_j\} = 1 - \alpha$ $(1 \leq j \leq k)$ when all $\mu_i = 0$, $i = 1, \ldots, k$.

This choice of constants enables control of the strong FWE. The step-up procedure is credited to Hochberg (1988).

**4. Properties of single-step.** Recall $b > 0$ and for $\Sigma$ intraclass, $-1/(k-1) < \rho$.

THEOREM 4.1. *For problems* (2.2) *and* (2.3), *suppose $\Sigma$ is intraclass. Suppose the loss function is* (2.7). *Then the single-step procedure is admissible if $\rho \geq -1/b$.*

For the proof see Appendix A.1.

REMARK 4.2. The proof that the single-step procedure is admissible under the given conditions is accomplished by demonstrating it is uniquely locally admissible in some sense. In a one-dimensional, one-sided hypothesis testing problem Lehmann (1986) describes a unique locally best test as one whose derivative of the power function evaluated at the null point is largest among all size $\alpha$ tests. For the multivariate global testing problem (2.1) a test is uniquely locally best in a direction if it has a similar property as in the one-dimensional case. Marden (1982) utilizes the notion of local admissibility for global testing problems. For our finite action problem we introduce a notion of unique local admissibility and demonstrate that the single-step procedure has this property. As in the global testing problem the focus is on the point $\boldsymbol{\mu} = \mathbf{0}$ and is linked to the derivative of a function of the risk evaluated at $\boldsymbol{\mu} = \mathbf{0}$.

Note that the theorem applies to the tree order model since in that case $\rho > 0$.

The admissibility result is particularly interesting in light of a result in CS (2004). There it is stated that when $\rho > 0$, no nontrivial Bayes test can be type-I monotone for problem (2.2). See CS (2004) for the definition of type-I monotone. The single-step procedure is type-I monotone, so a first guess might be that it is inadmissible. The result indicates that the first guess is incorrect.

Next we have:

THEOREM 4.3. *For problems* (2.2) *and* (2.3), *suppose $\Sigma = I$ (independence case). Suppose the loss function is* (2.7). *Then the single-step procedure is proper Bayes.*

For the proof see Appendix A.1.

This result should also be contrasted with the result in CS (2004) which states that the single-step procedure cannot be Bayes if $\Sigma$ is intraclass with $\rho > 0$:

Our final result of this section is:

THEOREM 4.4. *For problems* (2.2) *and* (2.3), *suppose $\Sigma$ is intraclass. Suppose the risk function for each component problem is* (2.10). *Then the single-step procedure is admissible for the vector risk VRI described in Result* 1.3.

For the proof see Appendix A.1.

Although the single-step procedure has the above desirable properties, many feel that single-step procedures are too conservative. That is, they do not detect significant effects often enough while controlling the FWE. Single-step procedures are somewhat akin to some simultaneous confidence bound procedures which are highly conservative, making it difficult to declare significance for an individual endpoint.

**5. Properties of step-down and step-up.** In this section the loss function is (2.7).

We establish the following theorems:

THEOREM 5.1. *For problems* (2.2) *and* (2.3) *for $\Sigma = I$, there exists a sequence of prior distributions for which the step-down procedure is a limit of a sequence of Bayes procedures.*

REMARK 5.2. In Theorem 4.1 it is shown that the single-step procedure has a limiting "local" optimality property. The limit point is **0**. In Theorem 5.1, however, it is shown that the step-down procedure has a limiting optimality property, but now the limiting parameter points receiving weight tend to $\infty$.

THEOREM 5.3. *For problem* (2.3), *$\Sigma = I$, $k = 2$, $b = 1$, the step-down procedure is admissible.*

THEOREM 5.4. *For problems* (2.2) *and* (2.3) *suppose $\Sigma$ is intraclass. Then the step-down procedure is admissible for vector risk VRI if and only if $\rho \geq 0$.*

For the proof see Appendix A.2.

Theorems 5.3 and 5.4 and the proofs of these theorems also apply to the step-up procedure given in Procedure 3.4. The most interesting properties for step-up are given in the companion paper CS (2005).

# APPENDIX

**A.1. Proofs of Theorems 4.1, 4.3 and 4.4.** In order to prove Theorem 4.1, we need a definition and theorem. First the definition.

For each $\mathbf{v} \in \Gamma$, let $\bar{\Omega}_\mathbf{v}$ be the closure of $\Omega_\mathbf{v}$. Let $R_\mathbf{v}(\boldsymbol{\psi}, \boldsymbol{\mu})$, for $\boldsymbol{\mu} \in \Omega_\mathbf{v}$, be the continuous extension of $R(\boldsymbol{\psi}, \boldsymbol{\mu})$ for $\boldsymbol{\mu} \in \bar{\Omega}_\mathbf{v}$. Note that the point $\mathbf{0} \in \bar{\Omega}_\mathbf{v}$ for all $\mathbf{v} \in \Gamma$. Since the risk function is continuous on each $\Omega_\mathbf{v}$, it follows that if $\boldsymbol{\psi}^*$ is better than $\boldsymbol{\psi}$, then $R_\mathbf{v}(\boldsymbol{\psi}^*, \boldsymbol{\mu}) \leq R_\mathbf{v}(\boldsymbol{\psi}, \boldsymbol{\mu})$ for all $\boldsymbol{\mu} \in \bar{\Omega}_\mathbf{v}$. In particular, $\boldsymbol{\psi}^*$ better than $\boldsymbol{\psi}$ implies $R_\mathbf{v}(\boldsymbol{\psi}^*, \mathbf{0}) \leq R_\mathbf{v}(\boldsymbol{\psi}, \mathbf{0})$ for all $\mathbf{v} \in \Gamma$.

The next theorem is useful when comparing decision procedures under the assumptions of this paper. That is, assume normality and assume the risk function is (2.11). In the case $k=1$ (which is the usual one-sided hypothesis testing problem), the theorem reduces to the well-known result that if $\boldsymbol{\psi}^*$ is better than $\boldsymbol{\psi}$, then their risks at zero (size of the test) must match.

THEOREM A.1. *If $\boldsymbol{\psi}^*$ is better than $\boldsymbol{\psi}$, then $R_\mathbf{v}(\boldsymbol{\psi}^*, \mathbf{0}) = R_\mathbf{v}(\boldsymbol{\psi}, \mathbf{0})$, all $\mathbf{v} \in \Gamma$.*

PROOF. The assumption that $\boldsymbol{\psi}^*$ is better than $\boldsymbol{\psi}$ implies that $R_\mathbf{v}(\boldsymbol{\psi}^*, \mathbf{0}) \leq R_\mathbf{v}(\boldsymbol{\psi}, \mathbf{0})$ for all $\mathbf{v} \in \Gamma$. Suppose $R_\mathbf{0}(\boldsymbol{\psi}^*, \mathbf{0}) < R_\mathbf{0}(\boldsymbol{\psi}, \mathbf{0})$. Since $R_\mathbf{0}(\boldsymbol{\psi}, \mathbf{0}) = E_\mathbf{0} \sum_{i=1}^k \psi_i(\mathbf{z})$ and $R_\mathbf{1}(\boldsymbol{\psi}, \mathbf{0}) = kb - E_\mathbf{0} \sum_{i=1}^k \psi_i(\mathbf{z})$, it follows that $R_\mathbf{1}(\boldsymbol{\psi}, \mathbf{0}) < R_\mathbf{1}(\boldsymbol{\psi}^*, \mathbf{0})$. This is a contradiction. A similar contradiction is reached if it is assumed that $R_\mathbf{1}(\boldsymbol{\psi}^*, \mathbf{0}) < R_\mathbf{1}(\boldsymbol{\psi}, \mathbf{0})$.

Now suppose for some $\mathbf{v} \in \Gamma_r$, $r = 1, \ldots, k-1$, $R_\mathbf{v}(\boldsymbol{\psi}^*, \mathbf{0}) < R_\mathbf{v}(\boldsymbol{\psi}, \mathbf{0})$, where $\Gamma_r = \{\mathbf{v} \in \Gamma : \sum_{i=1}^k v_i = r\}$. Then from (2.11)

$$(A.1) \qquad \sum_{\mathbf{v} \in \Gamma_r} R_\mathbf{v}(\boldsymbol{\psi}, \mathbf{0}) = \sum_{\mathbf{v} \in \Gamma_r} E_\mathbf{0}[\boldsymbol{\psi}'(\mathbf{1} - \mathbf{v}) + b(\mathbf{1} - \boldsymbol{\psi})'\mathbf{v}].$$

Now recognize that $\sum_{\mathbf{v} \in \Gamma_s} \mathbf{v} = \binom{k-1}{s-1}\mathbf{1}$ and collect terms so that (A.1) equals

$$(A.2) \qquad bk\binom{k-1}{r-1} + \left[\binom{k-1}{r} - b\binom{k-1}{k-r}\right] E_\mathbf{0} \sum_{i=1}^k \psi_i(\mathbf{z}).$$

If $[\binom{k-1}{r} - b\binom{k-1}{k-r}] > 0$, then $E_\mathbf{0} \sum_{v \in \Gamma_r} R_\mathbf{v}(\boldsymbol{\psi}^*, \mathbf{0}) < E_\mathbf{0} \sum_{\mathbf{v} \in \Gamma_r} R_\mathbf{v}(\boldsymbol{\psi}, \mathbf{0})$ implies $E_\mathbf{0} \sum_{i=1}^k \psi_i^*(z) < E_\mathbf{0} \sum_{i=1}^k \psi_i(\mathbf{z})$. This in turn implies that $R_\mathbf{1}(\boldsymbol{\psi}, \mathbf{0}) < R_\mathbf{1}(\boldsymbol{\psi}^*, \mathbf{0})$. This is a contradiction. If $[\binom{k-1}{r} - b\binom{k-1}{k-r}] < 0$ or equals 0, a similar contradiction is reached. Thus the theorem is proved. □

To prove Theorem 4.1 we need to study the behavior of linear combinations of the $R_\mathbf{v}$ functions. When $\Sigma$ is assumed to be intraclass we may write $\Sigma = \sigma^2((1-\rho)I + \rho \mathbf{1}\mathbf{1}')$. In this case

$$\Sigma^{-1} = (\sigma^2(1-\rho))^{-1}(I - G\mathbf{1}\mathbf{1})',$$

where

$$G = \rho/(1 + (k-1)\rho).$$

As earlier we take $\sigma^2 = 1$ without loss of generality.

PROOF OF THEOREM 4.1. Let $\boldsymbol{\psi}^*$ be the single-step procedure. Suppose $\boldsymbol{\psi}$ is better than $\boldsymbol{\psi}^*$. Then Theorem A.1 implies $\boldsymbol{\psi}$ cannot be uniformly better than the single-step procedure at $\mathbf{0}$; that is, there does not exist a $\boldsymbol{\psi}$ such that $R_{\mathbf{v}}(\boldsymbol{\psi}, \mathbf{0}) \leq R_{\mathbf{v}}(\boldsymbol{\psi}^*, \mathbf{0})$ for all $\mathbf{v}$, with strict inequality for some $\mathbf{v}$.

Therefore we need only consider procedures $\boldsymbol{\psi}$ such that

(A.3) $$R_{\mathbf{v}}(\boldsymbol{\psi}, \mathbf{0}) = R_{\mathbf{v}}(\boldsymbol{\psi}^*, \mathbf{0}) \qquad \text{for all } \mathbf{v} \in \Gamma.$$

For $\boldsymbol{\psi}$ satisfying (A.3) we study $\sum \lambda_{\mathbf{v}} R_{\mathbf{v}}(\boldsymbol{\psi}, \boldsymbol{\mu})$, where $\lambda_{\mathbf{v}}$ are coefficients that can depend on $\boldsymbol{\mu}$, for $\mathbf{v} \in \Gamma$ and where $R_{\mathbf{v}}(\boldsymbol{\psi}, \boldsymbol{\mu})$ is evaluated at $\boldsymbol{\mu} = \Delta \mathbf{v}$, $\Delta > 0$. In this case write $\sum_{\mathbf{v} \in \Gamma} \lambda_{\mathbf{v}} R_{\mathbf{v}}(\boldsymbol{\psi}, \boldsymbol{\mu}) = \sum_{\mathbf{v} \in \Gamma} \lambda_{\mathbf{v}}(\Delta) R_{\mathbf{v}}(\boldsymbol{\psi}, \Delta \mathbf{v})$. Among $\boldsymbol{\psi}$ satisfying (A.3), we show that $\boldsymbol{\psi}^*$ is the unique procedure that minimizes the derivative with respect to $\Delta$ of $\sum_{\mathbf{v} \in \Gamma} \lambda_{\mathbf{v}}(\Delta) R_{\mathbf{v}}(\boldsymbol{\psi}, \Delta \mathbf{v})$, evaluated at $\mathbf{0}$. This demonstrates the admissibility of $\boldsymbol{\psi}^*$.

Now we consider $\sum_{\mathbf{v} \in \Gamma} \lambda_{\mathbf{v}}(\Delta) R_{\mathbf{v}}(\boldsymbol{\psi}, \Delta \mathbf{v})$, which using (2.11) is

(A.4) $$\int \cdots \int \left\{ \sum_{s=0}^{k} \sum_{\mathbf{v} \in \Gamma_s} \lambda_{\mathbf{v}}(\Delta) [\mathbf{1}'(\boldsymbol{\psi} + b\mathbf{v}) - (1+b)\boldsymbol{\psi}'\mathbf{v}] f(\mathbf{z}|\mathbf{v}\Delta) \right\} d\mathbf{z},$$

where $f(\mathbf{z}|\mathbf{v}\Delta)$ is obtainable from

$$f(\mathbf{z}|\boldsymbol{\mu}) = (1/(2\pi)^{k/2} |\Sigma|^{1/2}) e^{-(1/2)(\mathbf{z}-\boldsymbol{\mu})' \Sigma^{-1} (\mathbf{z}-\boldsymbol{\mu})}.$$

For a chosen set of $\lambda_{\mathbf{v}}(\Delta)$ we seek a $\boldsymbol{\psi}$, among the class of procedures satisfying (A.3) that minimizes the derivative of (A.4) with respect to $\Delta$, evaluated at $\Delta = 0$.

Recall with $\sigma^2 = 1$, $\Sigma^{-1} = (1-\rho)^{-1}(I - G\mathbf{1}\mathbf{1}')$.

Now we choose $\lambda_{\mathbf{v}}$, $\mathbf{v} \in \Gamma$, as follows:

Let $\mathbf{C} = (C_1, \ldots, C_k)'$, where $C_i$ is given in Procedure 3.1.

Let $\boldsymbol{\varepsilon}_i = (0, \ldots, 1, 0, \ldots, 0)'$, that is, a vector with all zeros except 1 in the $i$th position. Let $\gamma = [(1 - Gk) + (1 + b)G]/b(1 - Gk) = (1 + b\rho)/b(1 - \rho)$, and note that $\gamma > 0$ if and only if $(1 + b\rho) > 0$. Let $\lambda_{\mathbf{0}}(\Delta) = 1$, $\lambda_{\mathbf{1}}(\Delta) = \gamma e^{-\mathbf{C}' \Sigma^{-1} \mathbf{1} \Delta}$, $\lambda_{\mathbf{v}}(\Delta) = e^{-\mathbf{C}' \Sigma^{-1} \mathbf{v} \Delta}$, for $\mathbf{v} = \boldsymbol{\varepsilon}_i$, $i = 1, \ldots, k$, and $\lambda_{\mathbf{v}}(\Delta) = 0$ otherwise.

The derivative of (A.4) with respect to $\Delta$ evaluated at $\Delta = 0$ is expressed as

$$\int \cdots \int \left\{ \sum_{s=1}^{k} \sum_{\mathbf{v} \in \Gamma_s} d(\lambda_{\mathbf{v}}(\Delta) f(\mathbf{z}|\mathbf{v}\Delta))/d\Delta|_{\Delta=0} \right.$$
$$\left. \times [\mathbf{1}'(\boldsymbol{\psi} + b\mathbf{v}) - (1+b)\boldsymbol{\psi}'\mathbf{v}] \right\} d\mathbf{z}$$

(A.5)
$$= (1/(1-\rho)) \int \cdots \int \left\{ \sum_{i=1}^{k} [\mathbf{1}'(\boldsymbol{\psi} + b\boldsymbol{\varepsilon}_i) - (1+b)\boldsymbol{\psi}'\boldsymbol{\varepsilon}_i] \right.$$
$$\times \boldsymbol{\varepsilon}_i'(I - G\mathbf{1}\mathbf{1}')(\mathbf{z} - \mathbf{C})$$
$$\left. + \gamma b(k - \boldsymbol{\psi}'\mathbf{1})k(1 - Gk)(\bar{z} - \bar{C}) \right\} f(\mathbf{z}|\mathbf{0}) \, d\mathbf{z}.$$

We will choose $\boldsymbol{\psi}(\mathbf{z})$ to minimize the integrand in (A.5) for each $\mathbf{z}$. Toward this end we evaluate the bracketed term on the right-hand side of (A.5), which becomes

(A.6)
$$\begin{aligned}
&\boldsymbol{\psi}'\mathbf{1}[k(\bar{z} - \bar{C}) - k^2 G(\bar{z} - \bar{C})] + b[k(\bar{z} - \bar{C}) - k^2 G(\bar{z} - \bar{C})] \\
&\quad - (1+b)\boldsymbol{\psi}'(\mathbf{z} - \mathbf{C}) + (1+b)kG\boldsymbol{\psi}'\mathbf{1}(\bar{z} - \bar{C}) \\
&\quad + \gamma b(k - \boldsymbol{\psi}'\mathbf{1})k(1 - Gk)(\bar{z} - \bar{C}) \\
&= -(1+b)\boldsymbol{\psi}'(\mathbf{z} - \mathbf{C}) \\
&\quad + \boldsymbol{\psi}'\mathbf{1}k(\bar{z} - \bar{C})\{(1 - Gk) + (1+b)G - \gamma b(1 - Gk)\} \\
&\quad + bk(\bar{z} - \bar{C})[(1 - G) + \gamma k(1 - Gk)].
\end{aligned}$$

At this point we recognize that by substituting the selected value of $\gamma$ in the bracketed term on the right-hand side of (A.6), the term becomes 0. Hence to minimize (A.6) we choose $\psi_i(\mathbf{z}) = 1$ if $z_i > C_i$ and choose $\psi_i(\mathbf{z}) = 0$ if $z_i < C_i$. This is the single-step procedure. Thus the single-step procedure is admissible for problem (2.2) [and for problem (2.3), since the same proof applies] if $\gamma > 0$. But $\gamma > 0$ if $(1 + b\rho) > 0$ which amounts to the given part of the theorem. □

In order to prove Theorem 4.3 we will need the following definition and theorem.

A decision procedure $\boldsymbol{\psi}^*$ is Bayes with respect to a prior distribution $\xi(\boldsymbol{\mu})$ if

$$E_\xi R(\boldsymbol{\psi}^*, \boldsymbol{\mu}) = \inf_{\boldsymbol{\psi}} E_\xi R(\boldsymbol{\psi}, \boldsymbol{\mu}).$$

In connection with Bayes procedures, let $q(\omega|\mathbf{z})$ denote the posterior probability of the subset $\omega \in \Omega$, given $\mathbf{z}$. Then the following theorem describes a Bayes procedure.

THEOREM A.2. *Consider the risk function in* (2.9). *The Bayes procedure is* $\boldsymbol{\psi}^* = (\psi_1^*, \ldots, \psi_k^*)'$, *where*

$$\psi_i^* = \begin{cases} 1, & \text{if } q(\Omega^{(i)}|\mathbf{z}) < b/(b+1), \\ 0, & \text{otherwise.} \end{cases}$$

PROOF. Since the loss function is additive, the sum of expected risks for the individual components is minimized by minimizing the expected risk for the individual components. The theorem follows by the same argument used for a single testing problem. See, for example, Mood, Graybill and Boes [(1974), page 417]. □

PROOF OF THEOREM 4.3. Choose a prior distribution such that $\mu_1, \ldots, \mu_k$ are independent. Then $q(\Omega^{(i)}|\mathbf{z})$ depends only on $z_i$. Furthermore, $q(\Omega^{(i)}|\mathbf{z})$ is a decreasing function of $z_i$. Use Theorem A.2 and the fact that the prior can be chosen so that $q(\Omega^{(i)}|\mathbf{z}) < b/(b+1)$ is equivalent to $Z_i > C_i$. □

PROOF OF THEOREM 4.4. We need only show that the procedure is componentwise admissible. That is, we need only prove that the test for each $H_i : \mu_i = 0$ vs $K_i : \mu_i > 0$ is admissible. It suffices to show that $Z_1 > C_1$ is an admissible test for $H_1 : \mu_1 = 0$ vs $K_1 : \mu_1 > 0$. To prove this we note that the multivariate normal density is proportional to

$$e^{-(1/2)\mathbf{z}'\Sigma^{-1}\mathbf{z}} e^{-(1/2)\boldsymbol{\mu}'\Sigma^{-1}\boldsymbol{\mu}} e^{\boldsymbol{\mu}'\Sigma^{-1}\mathbf{z}}. \tag{A.7}$$

Letting $\mathbf{y} = \Sigma^{-1}\mathbf{z}$, (A.7) can be written in exponential family form as

$$h(\mathbf{y})\beta(\boldsymbol{\mu})e^{\boldsymbol{\mu}'\mathbf{y}} = h(\mathbf{y})\beta(\boldsymbol{\mu})e^{y_1\mu_1 + \sum_{i=2}^k y_i\mu_i}.$$

A result of Matthes and Truax (1967) implies that any test of $H_1 : \mu_1 = 0$ vs $K_1 : \mu_1 > 0$ which is monotone in $y_1$ for fixed $(y_2, \ldots, y_k)$ is admissible. Here monotone means if $y_1' \leq y_1''$ and the test rejects for $y_1'$, then it must also reject for $y_1''$ when $y_2, \ldots, y_k$ are fixed.

Now note that the single-step procedure is of the form reject if $z_1 > C_1$. Since $\mathbf{z} = \Sigma\mathbf{y}$, this can be expressed as

$$y_1 + \rho \sum_{j=2}^k y_j > C_1. \tag{A.8}$$

From (A.8) we see that the test for $H_1$ is monotone in $y_1$ for fixed $(y_2, \ldots, y_k)$. □

**A.2. Proofs of Theorems 5.1, 5.3 and 5.4.**

PROOF OF THEOREM 5.1. A sequence of prior distributions will be put on various points of $\Omega_{\mathbf{v}}$. The amount of prior probability on each point will be expressed as a ratio where the denominator is always expressed as $D$ and $D$ is the sum of numerator terms. The sequence of priors is as follows: On $\Omega_{(0,\ldots,0)}$ the numerator of the prior probability is 1. On $\Omega_{(1,0,\ldots,0)}$ the numerator, $e^{\boldsymbol{\mu}'\boldsymbol{\mu}/2}e^{-C_k n^k}$, is put on $\mu_1 = n^k$; all other $\mu$'s are zero. On $\Omega_{(0,0,\ldots,0,1,0,\ldots,0)}$, where 1 is in the $i$th position, the numerator $e^{\boldsymbol{\mu}'\boldsymbol{\mu}/2}e^{-C_k n^k}$ is put on $\mu_i = n^k$; all other $\mu$'s are zero. On $\Omega_{(0,\ldots,0,1,0,\ldots,0,1,0,\ldots,0)}$, where 1 is in the $i$th and $j$th positions, the numerator $(1/2)e^{\boldsymbol{\mu}'\boldsymbol{\mu}/2}e^{-C_k n^k - C_{k-1} n^{k-1}}$ is put on the points $\mu_i = n^k$, $\mu_j = n^{k-1}$ and $\mu_i = n^{k-1}$, $\mu_j = n^k$ (all other $\mu$'s zero). On $\Omega_{(0,\ldots,0,1,0,\ldots,0,1,0,\ldots,0,1,0,\ldots,0)}$, where 1 is in the $i$th, $j$th and $\ell$th positions, the numerator $(1/3!)e^{\boldsymbol{\mu}'\boldsymbol{\mu}/2}e^{-C_k n^k - C_{k-1} n^{k-1} - C_{k-2} n^{k-2}}$ is put on six points, namely, $(\mu_i = n^k, \mu_j = n^{k-1}, \mu_\ell = n^{k-2})$, $(\mu_i = n^k, \mu_j = n^{k-2}, \mu_\ell = n^{k-1})$, $(\mu_i = n^{k-1}, \mu_j = n^k, \mu_\ell = n^{k-2})$, $(\mu_i = n^{k-1}, \mu_j = n^{k-2}, \mu_\ell = n^k)$, $(\mu_i = n^{k-2}, \mu_j = n^k, \mu_\ell = n^{k-1})$, $(\mu_i = n^{k-2}, \mu_j = n^{k-1}, \mu_\ell = n^k)$ (all other $\mu$'s are zero). In general, if $\mathbf{v} \in \Gamma_s$, then the numerator $(1/s!)e^{\boldsymbol{\mu}'\boldsymbol{\mu}/2}e^{\sum_{i=1}^{s} C_{k+1-i} n^{(k+1)-i}}$ is put on $s!$ points where the $\mu$'s are zero except for $(\mu_{j_1}, \ldots, \mu_{j_s})$ and all permutations where $\mu_{j_1}, \ldots, \mu_{j_s}$ correspond to $v_{j_1}, \ldots, v_{j_s}$ which are 1.

Next we indicate the numerators of posterior probabilities for each $\Omega_{\mathbf{v}}$. All posterior probabilities have the same denominator. We will note that for each fixed $\mathbf{z}$ one of the posterior probabilities will tend to 1. This fact means that the posterior risk will be minimized by choosing the action that corresponds to the $\Omega_{\mathbf{v}}$ whose posterior probability tends to 1. We will see that such a choice will correspond to the step-down procedure. Here are the numerators of the posterior probabilities denoted by $\xi(\Omega_{\mathbf{v}}|\mathbf{z})$. All denominators of the posterior probabilities are the same and the denominator is the sum of the numerators:

$$\xi(\Omega_{(0,\ldots,0)}|\mathbf{z}) = 1,$$

$$\xi(\Omega_{(1,0,\ldots,0)}|\mathbf{z}) = e^{(z_1 - C_k)n^k},$$

$$\xi(\Omega_{(0,\ldots,0,1,0,\ldots,0)}|\mathbf{z}) = e^{(z_i - C_k)n^k},$$

$$\xi(\Omega_{(0,\ldots,0,1,0,\ldots,0,1,0,\ldots,0)}|\mathbf{z}) = (1/2)[e^{(z_i - C_k)n^k + (z_j - C_{k-1})n^{k-1}} + e^{(z_j - C_k)n^k + (z_i - C_{k-1})n^{k-1}}].$$

For an arbitrary $\mathbf{v} \in \Gamma_s$,

$$(A.9) \quad \xi(\Omega_{\mathbf{v}}|\mathbf{z}) = (1/s!) \sum_{\substack{\text{all permutations} \\ \text{of } v_{j_1},\ldots,v_{j_s}}} \exp\left(\sum_{i=1}^{s}(z_{j_\ell} - C_{k+1-i})n^{(k+1)-i}\right),$$

where the indices $j_\ell$ reflect a permutation of $v_{j_1}, \ldots, v_{j_s}$.

At this point fix $\mathbf{z}$. Say, for example, and without loss of generality, $z_i > C_{k+1-i}$, $i = 1, 2, \ldots, r$, and $z_i < C_{k+1-i}$, $i = r+1, \ldots, k$. Then if $r = 0$, the posterior probability of $\Omega_{(0,\ldots,0)}$ denoted by $q(\Omega_{(0,\ldots,0)}|\mathbf{z})$ tends to 1 as $n \to \infty$. If $r \geq 1$, then $q(\Omega_{(1,\ldots,1,0,\ldots,0)}|\mathbf{z})$, with $r$ ones in $(1, \ldots, 1, 0, \ldots, 0)$, tends to 1 as $n \to \infty$. This is true since $\xi(\Omega_\mathbf{v}|\mathbf{z})$ tends to $\infty$ (except for $\Omega_{(0,\ldots,0)}$) as $n \to \infty$, but the ratio of $\xi(\Omega_{(1,\ldots,1,0,\ldots,0)}|\mathbf{z})/\xi(\Omega_\mathbf{v}|\mathbf{z})$ where $\mathbf{v}$ differs from $(1, \ldots, 1, 0, \ldots, 0)$ tends to $\infty$ as $n \to \infty$. Thus we have demonstrated that the step-down procedure is a limit of a sequence of Bayes procedures. □

PROOF OF THEOREM 5.3. For problem (2.3) the risk is taken from (2.10) and (2.11) except now $R_i(\psi_i, \boldsymbol{\mu}) = E_{\boldsymbol{\mu}}(\psi_i(\mathbf{z}))$ when $\mu_i \leq 0$. Let the step-down procedure be denoted by $\boldsymbol{\psi}^{\mathrm{SD}}$. Note that the risk function for an arbitrary procedure $\boldsymbol{\psi}$ is as follows:

For $\mu_1 > 0$, $\mu_2 > 0$,

$$R(\boldsymbol{\psi}, \boldsymbol{\mu}) = 2 - E_{\boldsymbol{\mu}}(\psi_1(\mathbf{z}) + \psi_2(\mathbf{z})). \tag{A.10}$$

For $\mu_1 \leq 0$, $\mu_2 \leq 0$,

$$R(\boldsymbol{\psi}, \boldsymbol{\mu}) = E_{\boldsymbol{\mu}}(\psi_1(\mathbf{z}) + \psi_2(\mathbf{z})). \tag{A.11}$$

For $\mu_1 > 0$, $\mu_2 \leq 0$,

$$R(\boldsymbol{\psi}, \boldsymbol{\mu}) = 1 - E_{\boldsymbol{\mu}}\psi_1(\mathbf{z}) + E_{\boldsymbol{\mu}}\psi_2(\mathbf{z}). \tag{A.12}$$

For $\mu_1 \leq 0$, $\mu_2 > 0$,

$$R(\boldsymbol{\psi}, \boldsymbol{\mu}) = 1 - E_{\boldsymbol{\mu}}\psi_2(\mathbf{z}) + E_{\boldsymbol{\mu}}\psi_1(\mathbf{z}). \tag{A.13}$$

Now if $\boldsymbol{\psi}^{\mathrm{SD}}$ is inadmissible from (A.10)–(A.13), then there exists a $\boldsymbol{\psi}^*$ which is better, that is,

$$E_{\boldsymbol{\mu}}(\psi_1^*(z) + \psi_2^*(z)) \geq E_{\boldsymbol{\mu}}(\psi_1^{\mathrm{SD}}(\mathbf{z}) + \psi_2^{\mathrm{SD}}(\mathbf{z})), \qquad \mu_1 > 0, \mu_2 > 0, \tag{A.14}$$

$$E_{\boldsymbol{\mu}}(\psi_1^{\mathrm{SD}}(z) + \psi_2^{\mathrm{SD}}(z)) \geq E_{\boldsymbol{\mu}}(\psi_1^*(\mathbf{z}) + \psi_2^*(\mathbf{z})), \qquad \mu_1 \leq 0, \mu_2 \leq 0, \tag{A.15}$$

$$E_{\boldsymbol{\mu}}(\psi_2^{\mathrm{SD}}(\mathbf{z}) - \psi_1^{\mathrm{SD}}(\mathbf{z})) \geq E_{\boldsymbol{\mu}}(\psi_2^*(z) - \psi_1^*(\mathbf{z})), \qquad \mu_1 > 0, \mu_2 \leq 0, \tag{A.16}$$

$$E_{\boldsymbol{\mu}}(\psi_1^{\mathrm{SD}}(\mathbf{z}) - \psi_2^{\mathrm{SD}}(\mathbf{z})) \geq E_{\boldsymbol{\mu}}(\psi_1^*(\mathbf{z}) - \psi_2^*(\mathbf{z})), \qquad \mu_1 \leq 0, \mu_2 > 0, \tag{A.17}$$

with at least one strict inequality for some $\boldsymbol{\mu}$. By letting either $\mu_1 \to 0$ or $\mu_2 \to 0$ or both $\mu_1 \to 0$, $\mu_2 \to 0$ in (A.14)–(A.17) we have that (A.14)–(A.17) hold whenever $\mu_1 \geq 0$, $\mu_2 \geq 0$; $\mu_1 \leq 0$, $\mu_2 \leq 0$; $\mu_1 \geq 0$, $\mu_2 \leq 0$; $\mu_1 \leq 0$, $\mu_2 \geq 0$, respectively.

Consider parameter points of the form $\boldsymbol{\mu} = (\mu_1, 0)'$, $\mu_1 \geq 0$. In this case (A.14) and (A.16) hold. Adding these two inequalities yields

$$E_{\boldsymbol{\mu}}\psi_1^*(\mathbf{z}) \geq E_{\boldsymbol{\mu}}\psi_1^{\mathrm{SD}}(\mathbf{z}). \tag{A.18}$$

Let $W(\mathbf{z}) = \psi_1^*(\mathbf{z}) - \psi_1^{\mathrm{SD}}(\mathbf{z})$ and let $\phi(u)$ be the standard normal density. Then (A.18) is

$$\begin{aligned}(\mathrm{A.19})\quad 0 &\leq \int_{-\infty}^\infty \int_{-\infty}^\infty W(\mathbf{z})\phi(z_1-\mu_1)\phi(z_2)\,dz_2\,dz_1 \\ &= \int_{-\infty}^\infty \int_{-\infty}^\infty W(\mathbf{z})\phi(z_2)\phi(z_1)e^{z_1\mu_1}e^{-C_2\mu_1}e^{C_2\mu_1}e^{-\mu_1^2/2}\,dz_2\,dz_1.\end{aligned}$$

Equivalently we have, for all $\mu_1 \geq 0$,

$$\begin{aligned}(\mathrm{A.20})\quad 0 &\leq \int_{-\infty}^\infty \int_{-\infty}^\infty W(\mathbf{z})\phi(z_2)\phi(z_1)e^{(z_1-C_2)\mu_1}\,dz_2\,dz_1 \\ &\leq 1 + \int_{C_2}^\infty \int_{-\infty}^\infty W(\mathbf{z})\phi(z_1)\phi(z_2)e^{(z_1-C_2)\mu_1}\,dz_2\,dz_1.\end{aligned}$$

In the last integral of (A.20) when $z_1 > C_2$, $\psi_1^{\mathrm{SD}}(\mathbf{z}) = 1$. Thus for any $\psi_1^*(\mathbf{z})$ $W(\mathbf{z}) \leq 0$ for all $z_1 \geq C_2$. If $W(\mathbf{z}) < 0$ on a set of positive Lebesgue measure, then the last integral in (A.20) tends to $-\infty$ as $\mu_1 \to \infty$. This would be a contradiction and so $W(\mathbf{z}) = 0$ a.s. for $z_1 > C_2$. Thus for $z_1 > C_2$, $\psi_1^*(\mathbf{z}) = \psi_1^{\mathrm{SD}}(\mathbf{z})$. This type of argument, letting $\mu \to \infty$ so that (A.20) $\to \infty$, is due to Stein. See, for example, Stein (1956).

In a similar fashion we show that

$$\begin{aligned}\psi_2^*(\mathbf{z}) &= \psi_2^{\mathrm{SD}}(\mathbf{z}) &&\text{for } z_2 > C_2, \\ \psi_1^*(\mathbf{z}) &= \psi_1^{\mathrm{SD}}(\mathbf{z}) &&\text{for } z_1 < C_1, \\ \psi_2^*(\mathbf{z}) &= \psi_2^{\mathrm{SD}}(\mathbf{z}) &&\text{for } z_2 < C_2.\end{aligned}$$

Therefore $\boldsymbol{\psi}^*(\mathbf{z}) = \boldsymbol{\psi}(\mathbf{z})$ whenever $(z_1 < C_1, z_2 < C_2)$, $(z_1 < C_1, z_2 > C_2)$, $(z_1 > C_2, z_2 < C_1)$ and $(z_1 > C_2, z_2 > C_2)$. That is, $\boldsymbol{\psi}^*(\mathbf{z}) = \boldsymbol{\psi}^{SD}(\mathbf{z})$ unless $C_1 < z_1 \leq C_2$ or $C_1 < z_2 < C_2$. Next return to (A.14) and consider $\boldsymbol{\mu} = (\mu_1, 1)'$. We have

$$(\mathrm{A.21})\quad 0 \leq E_{\boldsymbol{\mu}}\{\psi_1^*(\mathbf{z}) + \psi_2^*(\mathbf{z}) - \psi_1^{\mathrm{SD}}(\mathbf{z}) - \psi_2^{\mathrm{SD}}(\mathbf{z})\}.$$

Let $V(\mathbf{z}) = (\psi_1^*(\mathbf{z}) + \psi_2^*(\mathbf{z}) - \psi_1^{\mathrm{SD}}(\mathbf{z}) - \psi_2^{\mathrm{SD}}(\mathbf{z}))$. In the manner that (A.20) followed from (A.19), we have that (A.21) yields

$$\begin{aligned}(\mathrm{A.22})\quad 0 &\leq \int_{-\infty}^\infty \int_{-\infty}^\infty V(\mathbf{z})\phi(z_1)\phi(z_2-1)e^{(z_1-C_1)\mu_1}\,dz_2\,dz_1 \\ &\leq 2 + \int_{C_2}^\infty \int_{C_1}^{C_2} V(\mathbf{z})\phi(z_1)\phi(z_2-1)e^{(z_1-C_1)\mu_1}\,dz_2\,dz_1.\end{aligned}$$

When $z_1 > C_2$ and $z_2 \in [C_1, C_2]$ we have $V(\mathbf{z}) \leq 0$. As before, we have a contradiction in (A.22) as $\mu_1 \to \infty$ unless $\boldsymbol{\psi}^*(\mathbf{z}) = \mathbf{1}$ in $\{\mathbf{z}: z_1 > C_2, C_1 < z_2 < C_2\}$.

Similarly it can be shown that $\boldsymbol{\psi}^*(\mathbf{z}) = \boldsymbol{\psi}^{SD}(\mathbf{z})$ for almost all $\mathbf{z}$ not lying in the box

$$(\mathrm{A.23})\quad \{\mathbf{z}: C_1 \leq z_1 \leq C_2, C_1 \leq z_2 \leq C_2\}.$$

The final step is to show that $\boldsymbol{\psi}^*(\mathbf{z}) = \boldsymbol{\psi}^{\mathrm{SD}}(\mathbf{z})$ on (A.23). Now $\psi^{\mathrm{SD}}(\mathbf{z}) = 1$ on (A.23), so (A.15) would be violated when $\boldsymbol{\mu} = \mathbf{0}$ if $\boldsymbol{\psi}^*(\mathbf{z}) \neq \mathbf{1}$ on a set of positive measure in (A.23). This completes the proof. $\square$

PROOF OF THEOREM 5.4. As in the case of the proof of Theorem 4.4 we appeal to the Matthes and Truax (1967) theorem. We must show that the step-down procedure is monotone in $y_1$ for fixed $(y_2, \ldots, y_k)$ if and only if $\rho \geq 0$.

Now recognize that the step-down procedure is of the form reject $H_1$ if $z_1 > C(\mathbf{z}^{(2)})$, $\mathbf{z}^{(2)} = (z_2, \ldots, z_k)'$, or in terms of the coordinates of $\mathbf{y}$ it is of the form reject if

$$
\begin{aligned}
(\mathrm{A.24}) \quad y_1 + \rho \sum_{i=2}^{k} y_i > C\bigg( & \rho y_1 + y_2 + \rho \sum_{i=3}^{k} y_i, \\
& \rho(y_1 + y_2) + y_3 + \rho \sum_{i=4}^{k} y_i, \ldots, \rho \sum_{i=1}^{k-1} y_i + y_k \bigg).
\end{aligned}
$$

Note that the left-hand side of (A.24) is increasing in $y_1$ for fixed $y_2, \ldots, y_k$. We claim the right-hand side of (A.24) is nonincreasing in $y_1$ for fixed $y_2, \ldots, y_k$ as long as $\rho \geq 0$. To see this, note that $C(\mathbf{z}^{(2)})$ is a nonincreasing function of its arguments. That is, as any $z_i$, $i = 2, \ldots, k$, increases it becomes easier to reject $H_1$; that is, the critical value in the step-down sequence can only remain the same or become smaller. For example, if all $z_i$, $i = 2, \ldots, k$, are less than $C_k$, then $C(\mathbf{z}^{(2)}) = C_k$. If exactly one of $z_2, \ldots, z_k$ is bigger than $C_k$, then $C(\mathbf{z}^{(2)}) = C_{k-1}$. Thus the conditions of the Matthes and Truax theorem are met and the step-down procedure (and step-up procedure) are admissible for VRI as long as $\rho \geq 0$.

Next we show that if $\rho < 0$, then the step-down (step-up) procedures are not monotone on some sections (monotone in $y_1$ for fixed $y_2, \ldots, y_k$) and therefore can be improved on these sections. Toward this end recall that

$$(\mathrm{A.25}) \qquad \mathbf{y} = \Sigma^{-1}\mathbf{z} \quad \text{and} \quad \mathbf{z} = \Sigma \mathbf{y}.$$

Let $\mathbf{r} = (1 \ \rho \ \rho \cdots \rho)'$ be the first column of $\Sigma$ and define the points $\mathbf{z}^*$ and $\mathbf{z}^{**}$ as follows:

$$\mathbf{z}^* = ((C_{k-1} + C_k)/2, \ C_k, \ldots, C_k),$$

$$\mathbf{z}^{**} = \mathbf{z}^* - \varepsilon \mathbf{r}.$$

The step-down procedure accepts $H_1$ when $\mathbf{z}^*$ is observed (it is a boundary point of the acceptance region). Since $\rho < 0$ when $0 < \varepsilon$ is sufficiently small, the step-down procedure will reject $\mathbf{z}^{**}$. It follows from (A.25) that

$$\begin{aligned}
\mathbf{y}^* = \Sigma \mathbf{z}^* &= \Sigma(\mathbf{z}^{**} + \varepsilon \mathbf{r}) = \Sigma \mathbf{z}^{**} + \varepsilon \Sigma \mathbf{r} \\
&= \mathbf{y}^{**} + (\varepsilon, 0, \ldots, 0)'.
\end{aligned}$$

Thus the step-down is not monotone in $y_1$. A similar argument works for step-up.

□

**Acknowledgment.** We are grateful to the Associate Editor who suggested a better version of Theorem 4.1 and who supplied a more useful weight function to be used in the proof.

Department of Statistics
Rutgers—The State University
 of New Jersey
Hill Center, Busch Campus
Piscataway, New Jersey 08854
USA
e-mail: artcohen@rci.rutgers.edu
e-mail: sackrowi@rci.rutgers.edu